\documentclass{amsart}
\usepackage[english]{babel}
\usepackage[utf8]{inputenc}
\usepackage{times}
\usepackage[T1]{fontenc}
\usepackage{amsmath,amsthm,amssymb}
\usepackage{xcolor}
\usepackage{lineno}
\usepackage{tikz}
\usetikzlibrary{decorations.pathmorphing,positioning,patterns}
\ifx\UseOption\undefined
\def\UseOption{pdftikz,dessinempost,color}
\fi
\usepackage{optional}
\opt{color}{
\colorlet{couleurlien}{green!20!black}
\colorlet{couleurlink}{red!60!black}
}
\opt{bw}{
\colorlet{couleurlien}{black}
\colorlet{couleurlink}{black}
}

\usepackage[pdftex,colorlinks,citecolor=couleurlien,linkcolor=couleurlink]{hyperref}

\newcommand\sshift[1]{\ensuremath{\mathbf{#1}}}

\newcommand\couleur[1]{\mathbf{#1}}
\newcommand\alphab{\ensuremath{\couleur{\Sigma}}}
\newcommand\motif[1]{\ensuremath{\mathbf{#1}}}

\newcommand\domain[1]{\ensuremath{\mathcal{D}_{#1}}}

\newcommand\espace[1][\plan]{\ensuremath{\alphab^{#1}}}
\newcommand\myspace{\espace[{\Z}^d]}
\newcommand\conf[1]{#1}

\newcommand\ForbPat{\ensuremath{\mathcal{F}}}

\newcommand\card[1]{\ensuremath{\lvert{}#1\lvert{}}}

\newcommand\layer[1]{\ensuremath{L_{#1}}}

\newcommand\latin[1]{\emph{#1}}

\newcommand\ie{\latin{i.e.}, }
\newcommand\eg{\latin{e.g.}, }
\newcommand\Z{\ensuremath{\mathbb{Z}}}
\newcommand\N{\ensuremath{\mathbb{N}}}

\newcommand\machine[1]{\ensuremath{\mathcal{M}_{#1}}}
\newcommand\massp[1]{\ensuremath{\mathcal{#1}}}
\newcommand\massred{\ensuremath{\leq}}
\newcommand\massredrev{\ensuremath{\geq}}
\newcommand\degdif[1]{\ensuremath{\mathbf{#1}}}
\newcommand\mindegdif{\ensuremath{\degdif{0}}}
\newcommand\maxdegdif{\ensuremath{\degdif{1}}}
\newcommand\masstrad{\ensuremath{\Theta}}
\newcommand\massdetrad{\ensuremath{\Theta^{-1}}}

\newcommand\codeconfmot[2]{\ensuremath{#1^{#2}}}
\newcommand\codemotconf[2]{\ensuremath{#1^{#2}}}
\newcommand\ucfunc{\ensuremath{\mathcal{U}}}
\newcommand\recop[1]{\ensuremath{\Psi_{#1}}}
\newcommand\effofcode[1]{\ensuremath{\sshift{X}_{#1}}}
\newcommand\masspbofset{\ensuremath{\mathcal{M}_{\textrm{set}}}}
\newcommand\masspbofsubshift{\ensuremath{\mathcal{M}_{\textrm{subshift}}}}
\newcommand\medvedevsup{\ensuremath{\vee}}
\newcommand\medvedevinf{\ensuremath{\wedge}}
\newcommand\entrop[1]{\ensuremath{h(#1)}}
\newcommand\univdef[3]{\ensuremath{(#1,#2,#3)-}universal}
\newcommand\univclas{\univdef{\mathcal{R}}{\mathcal{G}}{\mathcal{B}}}
\newcommand\univclasinf{\univdef{\mathcal{R}}{\mathcal{G}}{\infty}}

\newcommand\val[2]{\ensuremath{\conf{#1}_{#2}}}

\newcommand\mgreen{green!75!black!50}

\newcommand\leftbrack[2]{
\opt{bw}{
\filldraw[fill=black!20] (#1,#2) rectangle +(1,1);
}
\opt{color}{
\filldraw[fill=purple] (#1,#2) rectangle +(1,1);
}
\draw (#1,#2+.5) -- +(1,.5);
\draw (#1,#2+.5) -- +(1,-.5);
}

\newcommand\rightbrack[2]{
\opt{bw}{
\filldraw[fill=black!20] (#1,#2) rectangle +(1,1);
}
\opt{color}{
\filldraw[fill=purple] (#1,#2) rectangle +(1,1);
}
\draw (#1+1,#2+.5) -- +(-1,.5);
\draw (#1+1,#2+.5) -- +(-1,-.5);
}
\newcommand\codingcell[2]{
\opt{bw}{
\filldraw[fill=black!20] (#1,#2) rectangle +(1,1);
}
\opt{color}{
\filldraw[fill=\mgreen] (#1,#2) rectangle +(1,1);
}
}
\newcommand\obscell[2]{
\opt{bw}{
\filldraw[fill=black!40] (#1,#2) rectangle +(1,1);
}
\opt{color}{
\filldraw[fill=orange] (#1,#2) rectangle +(1,1);
}
}

\newcommand\bijnz[1]{\ensuremath{\mathcal{B}_{#1}}}

\newcommand\projsub[3]{\ensuremath{#1^{(#2,#3)}}}

\newtheorem{theorem}{Theorem}[section]
\newtheorem{corollary}{Corollary}

\newtheorem{lemma}[theorem]{Lemma}

\theoremstyle{definition}
\newtheorem{definition}[theorem]{Definition}

\newtheorem{op}{Problem}

\title[Universality and Medvedev degrees]{Universality in symbolic dynamics constrained by Medvedev degrees}
\author[Alexis Ballier]{}
\email{aballier@dim.uchile.cl}
\thanks{This research has been supported by the FONDECYT Postdoctorado Proyecto
3110088}


\begin{document}

\maketitle

\centerline{\scshape Alexis Ballier}
\medskip
{\footnotesize
        \centerline{Centro de Modelamiento Matemático}
        \centerline{Av. Blanco Encalada 2120 Piso 7, Oficina 710}
\centerline{Santiago, Chile}
} 

\bigskip


\begin{abstract}
        We define a weak notion of universality in symbolic
        dynamics and, by generalizing a proof of Mike Hochman, we prove that
        this yields necessary conditions on the forbidden
        patterns defining a universal subshift: These forbidden patterns are
        necessarily in a Medvedev degree greater or equal than the degree of the
        set of subshifts for which it is universal.
        
        We also show that this necessary condition is optimal by giving constructions of
        universal subshifts in the same Medvedev degree as the subshifts they
        simulate and prove that this universality can be achieved by the sofic
        projective subdynamics of such a subshift as soon as the necessary
        conditions are verified.

        This could be summarized as: There are obstructions for the
        existence of universal subshifts due to the theory of computability and
        they are the only ones.
\end{abstract}

An important problem in topological dynamics is to understand when one system
can ``simulate'' another one: For example, two systems will be essentially the
same if they are conjugate, a system will be simpler than another one if the
former is a factor of the latter, etc.
In symbolic dynamics, most of the one-dimensional cases have been settled by
algebraic invariants: An irreducible SFT factors onto another SFT of strictly
smaller entropy if and only if the necessary conditions on the periodic points are
met~\cite{boylelowentropfact}, and recent results have settled the sofic
case~\cite{kriegersof}. On the other hand, even if we know some sufficient
conditions~\cite{jmaddfactor}, the situation tends to be much more
complex in multi-dimensional symbolic dynamics~\cite{msmbrp}.

In theoretical computer science, a fundamental result is the existence of a
universal Turing machine~\cite{turingent}: This particular machine enables
constructing a good part of the theory of computability~\cite{rogersbook} and
will be used implicitly in this article.
While the one-dimensional SFTs are quite simple and well
understood~\cite{marcuslind}, multi-dimensional SFTs can exhibit computational
properties~\cite{bergerthesis,robinson}.
It then becomes natural to ask whether there exist multi-dimensional SFTs that
can ``simulate'' every other SFT, that is, as in the Turing machine case, a
universal SFT.
Unfortunately, Mike Hochman proved that there cannot exist such
an SFT~\cite{hochmanuniv}.

Usually, the notion of simulation between subshifts is the factoring
relation, or sofic projective subdynamics~\cite{msrpsubdyn} if we want more
flexibility.
In this article, we extend this notion of simulation to any computable
function that preserves the subshift structure: Computability between subshift
is modelled by \emph{recursive operators} and preserving the structure is what
we define as \emph{subshift-preserving}.
A subshift will then simulate another one when there exists a \emph{computable}
function (a recursive operator) mapping the former onto a \emph{subsystem} of
the latter.
Recursive operators allow to define a pre-order between sets of functions (a
configuration of a subshift may be seen as a function) known as the
\emph{Medvedev reduction} and its equivalence classes are the \emph{Medvedev
degrees}~\cite{medvedev} (or \emph{degrees of difficulty}~\cite{rogersbook}).
The idea behind Medvedev degrees is that a set will be \emph{simpler} than
another one if from any element of the latter we can compute an element of the
former. For example, if we consider factor maps, a subshift will be simpler
than another one if the latter factors onto a \emph{subsystem} of the former.
Based on these relations we define a notion of universality in symbolic
dynamics (Definition~\ref{defi:univ}): A subshift is universal for a set of
subshifts if it can \emph{simulate}, in that way, a subsystem of every subshift in
the class.
Being honest, this definition is purely academical and
is not what one would consider a natural definition of universality.
However, this notion includes most of what one would call
universality: A universal subshift could be a subshift that factors onto all the
subshifts of a given set, that admits all these subshifts as sofic projective
subdynamics or even that admits a rescaling factoring on a rescaling of every
subshift in the set such as in \cite{gtlbunivac}. All these notions fall under
our definition of universality.

Even with this weak notion of universality, we are able to generalize M.
Hochman's proof~\cite{hochmanuniv} to prove that the Medvedev degree of the
class of simulated subshift cannot be greater than the difficulty of defining
the universal subshift (Theorems~\ref{thm:univgtclass}
and~\ref{thm:univgtclasssft}).
We prove that these conditions are optimal by constructing universal subshifts
of the smallest possible Medvedev degree: If we allow any computable
function as a simulation then we can always construct a one-dimensional
universal subshift (Theorem~\ref{thm:univ1d}). If we restrict the simulation
to sofic projective subdynamics then the simulation cannot increase the
dimension and we are able to prove that universality can be achieved by sofic
projective subdynamics as soon as the necessary conditions are met
(Theorem~\ref{thm:existsunivdimqcqdown}).
If we only allow factor maps then entropy has to be bounded since it cannot
increase with factor maps but then we do not know if, additionally, bounding
the entropy is sufficient for the existence of such a universal subshift.

\section{Definitions, notations and codings}

\subsection{Elements of symbolic dynamics}

Let $\alphab$ be a finite set, called the \emph{alphabet}.
For an integer $d$, $\myspace$ is the set of configurations of dimension $d$
over $\alphab$, \ie the functions from ${\Z}^d$ to $\alphab$, called the
\emph{fullshift} of dimension $d$ over $\alphab$.
If $\conf{x}\in\myspace$ is a configuration, for $i\in\Z^d$, we denote by
$\val{x}{i}$ the value of $\conf{x}$ at $i$.
We endow $\alphab$ with the discrete topology and $\myspace$ with its product
which turns it in a compact metrizable space.
A pattern $\motif{P}$ is a mapping from a finite square subset
$\domain{P}=[-n;n]^d$ of $\Z^d$ to $\alphab$.
A pattern $\motif{P}$ is said to be \emph{supported} over $D$ if
$\domain{P}\subseteq{}D$.
For $D\subseteq\Z^d$, $x_{|D}$ denotes the restriction of $x:\Z^d\to\alphab$ to
$x_{|D}:D\to\alphab$.
The \emph{shift} of vector $i$, for $i\in\Z^d$, is defined by:
$\val{\sigma_i(\conf{x})}{j}=\val{x}{i+j}$, for every $j\in\Z^d$.
A subset $\sshift{X}$ of $\myspace$ is said to be \emph{shift-invariant} if for
all $i\in\Z^d$ and all $x\in\sshift{X}$, $\sigma_i(x)$ belongs to $\sshift{X}$.
A \emph{subshift} is a \emph{closed and shift-invariant} subset of $\myspace$.

A configuration $\conf{c}\in\myspace$ is said to contain a pattern $\motif{P}$
if there exists $i\in\Z^d$ such that
$\sigma_i(\conf{c})_{|\domain{P}}=\motif{P}$.
Given a set of patterns $\ForbPat$, $\sshift{X}_{\ForbPat}$ is the
\emph{subshift} defined by forbidding $\ForbPat$:
The set of configurations that contain no pattern of $\ForbPat$.
Such $\sshift{X}_{\ForbPat}$'s are closed and
shift-invariant; for any closed and shift-invariant subset of $\myspace$
there exists such an $\ForbPat$ defining it~\cite{Hedlund69}.
For a subshift $\sshift{X}$, if there exists a recursively enumerable $\ForbPat$
such that $\sshift{X}=\sshift{X}_{\ForbPat}$ then $\sshift{X}$
is an \emph{effective subshift}; if $\ForbPat$ is finite
then $\sshift{X}$ is a \emph{subshift of finite type} (\emph{SFT} in
short).

An onto function $\varphi:\sshift{X}\to\sshift{Y}$ between two subshifts is
called a \emph{factor map} if it is continuous and
shift commuting (\ie $\varphi\circ\sigma_i=\sigma_i\circ\varphi$ for all
$i\in\Z^d$).
By Hedlund's theorem~\cite{Hedlund69}, factor maps correspond to
sliding block codes: For 
$\varphi:\sshift{X}\to\sshift{Y}$ where $\sshift{X}\subseteq{}\myspace$ and
$\sshift{Y}\subseteq{}\couleur{\Gamma}^{\Z^d}$ 
there exists a finite subset $D$ of $\Z^d$ and a function
$\Phi:\alphab^{D}\to\couleur{\Gamma}$ such that for any $x\in\sshift{X}$,
$\varphi(x)_i=\Phi(\sigma_i(x)_{|D})$; conversely, a sliding block code is a factor map
between its domain and its image.

A subgroup $G$ of $\Z^d$ is isomorphic to some $\Z^k$ for $k\leq{}d$. Denote
this isomorphism $h:\Z^k\to{}G$.
Considering a subshift $\sshift{X}$ over $\myspace$, we can define its
\emph{projective subdynamics over $G$} and see it as a subshift of
$\alphab^{\Z^k}$ by $H$ defined by $H(x)_i=x_{h(i)}$.
Given, in addition, a factor map $\varphi$ with domain $\sshift{X}$,
we can define the \emph{$(h,\varphi{})$
sofic projective subdynamics} of $\sshift{X}$ by
$\projsub{\sshift{X}}{h}{\varphi}=H\circ\varphi(\sshift{X})$~\cite{hochmanrecsshift,msrpsubdyn}.
Remark that $\projsub{\sshift{X}}{h}{\varphi}$ can have greater entropy than
$\sshift{X}$ even if $G$ is a subgroup of finite index of $\Z^d$, and in this
case the operation is very similar to the bulking used for
simulations between cellular automata~\cite{bulkingi,bulkingii,gtlbunivac}.
The case where $G$ is a subgroup of finite index of $\Z^d$, meaning the
projective subdynamics do \emph{not} decrease the dimension, is maybe the most
natural notion of sofic projective subdynamics in the context of universality
since it avoids the entropy constraints while remaining close to a factor map.
When $G$ is a subgroup of finite index of $\Z^d$ we say that sofic
projective subdynamics over $G$ are \emph{finite index sofic projective
subdynamics}.

\subsection{Mass problems and subshifts}

In this paper, we want to study the \emph{Medvedev
degrees}~\cite{medvedev,rogersbook} of subshifts and the links that can be
obtained between universality and the relations among those degrees.
A \emph{mass problem} is a collection of total functions from $\N$ to $\N$, \ie
a subset of $\N^{\N}$.
We can assume, without loss of generality, that subshifts are defined on an
alphabet that is a subset of $\N$.
The collection of all configurations of any dimension can then be
effectively coded and decoded as elements of $\N^{\N}$: For an integer $d$, let
$\bijnz{d}$ be a (fixed) computable bijection from $\N$ onto $\Z^d$;
$\conf{c}\in\myspace$ is coded by $\codeconfmot{w}{\conf{c}}\in\N^{\N}$ such that:
$\val{\codeconfmot{w}{\conf{c}}}{0}=\card{\alphab},
\val{\codeconfmot{w}{\conf{c}}}{1}=d, \val{\codeconfmot{w}{\conf{c}}}{n+2} =
\val{\conf{c}}{\bijnz{d}(n)}$.
The decoding of $w\in\N^{\N}$ is done with the same idea:
Let $\alphab = \left\{0,\ldots,\val{w}{0}-1\right\}$, $d=\val{w}{1}$, and then
decode $w\in\N^{\N}$ as $\codemotconf{\conf{c}}{w}\in\myspace$ such that
$\val{\codemotconf{\conf{c}}{w}}{i}=\val{w}{\bijnz{d}^{-1}(i)+2}\mod\card{\alphab}$.
Denote by $\masstrad$ the function that codes $\conf{c}$ as
$\codeconfmot{w}{\conf{c}}\in\N^{\N}$ and $\massdetrad$ the function that
decodes $w\in\N^{\N}$ as $\codemotconf{\conf{c}}{w}$.
For any $\alphab$ and any $d$, $\masstrad$ is a bijection between $\myspace$ and
$\masstrad(\myspace)$ and we denote by $\massdetrad$ its inverse.
As such, any subshift can (effectively) be seen as a mass problem \latin{via}
the $\masstrad$ function.

Since coding and decoding between $\myspace$ and $\N^{\N}$ should really be
seen as a technical detail, we say that a mass problem $\mathcal{C}$ is a
subshift, SFT or effective subshift if $\massdetrad(\mathcal{C})$ is a subset of
some $\myspace$ that is, respectively, a subshift, SFT or effective subshift.

\subsection{Recursive operators}

\emph{Turing machines}~\cite{turingent} are often taken as the standard model
of computation. We mainly use the, more
abstract, recursive operators defined below. 
The reader who is not familiar with the concept of Turing machines may understand
it as an algorithm in any programming language, or refer to standard textbooks
on the subject such as \cite{rogersbook,odifreddi1,sbcoopercomp}.
While (classical) Turing machines deal with finite words on a finite alphabet,
with mass problems we have to
deal with infinite words over a countable alphabet (elements of $\N^{\N}$),
hence we define the analogous of Turing machines in this context:

\begin{definition}
A partial function $\Psi:\N^{\N}\to\N^{\N}$ is said to be a \emph{recursive
operator} if there exists an oracle Turing machine $\machine{\Psi}$ such that
for any $f\in\N^{\N}$ where $\Psi$ is defined,
for any $n\in\N$, $\machine{\Psi}$ given as oracle $f$ and input $n$ eventually
halts with output $\val{\Psi(f)}{n}\in\N$.

$\machine{\Psi}$ is said to \emph{compute $\Psi$}.
\end{definition}

Recursive operators are continuous for the natural
topology over $\N^{\N}$: The Turing machine $\machine{\Psi}$ computing $\Psi$,
given any input $n$, halts after a finite number of steps and thus needs only a
finite number of values from $f$; for any integer $n$, there
exists $k_n$ such that the output of the machine is independent of the values of
$\val{f}{i}$, for $i\geq{}k_n$.
Therefore, $\Psi$ is a continuous function.

Since Turing machines are recursively enumerable, (partial) recursive operators
are recursively enumerable too; in the remaining of this article, we fix such a
recursive enumeration and denote by
$\recop{n}$ the $n^{th}$ recursive operator in this fixed enumeration.

A factor map 
$\varphi{}:\couleur{\Sigma}^{\Z^d}\to\couleur{\Pi}^{\Z^d}$ is a computable
function since it is a sliding block code~\cite{Hedlund69}: See, \eg
\cite[Theorem~$1$]{simpsonmedv}.
$\psi=\masstrad\circ\varphi\circ\massdetrad:\masstrad(\myspace)\to\masstrad(\couleur{\Pi}^{\Z^d})$
can, therefore, be seen as a partial recursive operator.
By the same type of arguments,
the sofic projective subdynamics operation
($\sshift{X}\to\projsub{\sshift{X}}{h}{\varphi}$) is a partial recursive
operator.

\begin{definition}
        For a subshift $\sshift{X}$, a recursive operator $\Psi$ is said to be
        \emph{$\sshift{X}$-subshift-preserving} if $\Psi$ is defined on (the
        coding by $\masstrad$ of) $\sshift{X}$ and its image is a subshift.
\end{definition}

For example, factor maps, sofic projective subdynamics or the $\preceq$
simulation from \cite{gtlbunivac} give
subshift-preserving recursive operators when transformed as such.
Since recursive operators are functions from $\N^{\N}$ to $\N^{\N}$ and we will
often have to apply them to the coding of subshifts, when there is no ambiguity,
for a subshift $\sshift{X}$, we will denote $\recop{n}(\sshift{X})$ for
$\massdetrad\circ\recop{n}\circ\masstrad{}(\sshift{X})$.

\subsection{Medvedev degrees}

A mass problem $\massp{A}$ is said to be \emph{Medvedev reducible} to the mass
problem $\massp{B}$ if there exists a recursive operator $\Psi_n$ such that
$\Psi_n$ is defined over $\massp{B}$ and
$\Psi_n(\massp{B})\subseteq\massp{A}$.
We denote this by $\massp{A}\massred_n\massp{B}$. We write
$\massp{A}\massred\massp{B}$ if there exists $n$ such that
$\massp{A}\massred_n\massp{B}$.
Intuitively, $\massp{A}\massred\massp{B}$ means that there exists an automatic
procedure such that when given any solution of a problem $\massp{B}$ finds a
solution to the problem $\massp{A}$.
The relation $\massred$ is a pre-order. When quotiented by its induced
equivalence relation ($\massred\wedge\massredrev$), the set of mass problems
becomes the \emph{Medvedev lattice}~\cite{medvedev} (also known as the
\emph{degrees of difficulty}~\cite{rogersbook}).

The empty mass problem is in the highest degree of
difficulty, denoted by $\maxdegdif$; a mass problem containing a computable
function is always in the minimal degree, denoted by $\mindegdif$, and
any mass problem in $\mindegdif$ contains a computable function.
Simpson showed that the sublattice of the $\Pi_1^0$ subsets of
$\left\{0,1\right\}^{\N}$ is exactly the Medvedev lattice of $2-$dimensional
SFTs~\cite{simpsonmedv}.
Remark that if a subshift $\sshift{X}$ factors onto
$\sshift{Y}$ then $\sshift{Y}\massred\sshift{X}$ since
factor maps can be seen as recursive operators.
The Medvedev degree is therefore a conjugacy
invariant.

Medvedev degrees ordered by the relation $\massred$ form a
lattice: There is a $\sup$, denoted by $\medvedevsup$ and representing
the product of the two sets, and an $\inf$ denoted by $\medvedevinf$ and
representing the disjoint union of the two sets; these are defined formally by:
\[
\begin{array}{rcl}
\massp{A}\medvedevsup\massp{B}& = & \left\{x\in\N^{\N}| \exists a\in\massp{A},
\exists b\in\massp{B}, \forall n\in\N, x(2n)=a(n), x(2n+1)=b(n) \right\}\\
i\oplus\massp{A}&=&\left\{x\in\N^{\N}| x(0)=i, \exists a\in\massp{A}, \forall n\in\N, x(n+1)=a(n)\right\}\\
\massp{A}\medvedevinf\massp{B}& = & \left(0\oplus\massp{A}\right)\cup
\left(1\oplus\massp{B}\right)\\
\end{array}
\]
Strictly speaking,
$\left(\massp{A}\medvedevsup\massp{B}\right)\medvedevsup\massp{C}$ is different
from $\massp{A}\medvedevsup\left(\massp{B}\medvedevsup\massp{C}\right)$, but 
they are in the same Medvedev degree (the same applies for
$\medvedevinf$).
Therefore, we can use the $\medvedevsup$ and $\medvedevinf$ operations
associatively on Medvedev degrees.

\subsection{Codings}

We want to compare the Medvedev degrees of sets of subshifts and those of the
forbidden patterns defining them;
therefore we also need to define how
to code these sets as mass problems.

We can associate to any subset $\mathcal{N}$ of $\N$ a canonical mass problem
$\masspbofset(\mathcal{N})\subseteq\N^\N$ defined by:
$
\masspbofset(\mathcal{N}) = \left\{x\in\N^\N| n \in \mathcal{N}
\Leftrightarrow \exists i\in\N, \val{x}{i}=n \right\}
$.
$\mathcal{N}$ is recursively enumerable if and only if
$\masspbofset(\mathcal{N})$ contains a computable point, thus if and only if
$\masspbofset(\mathcal{N})$ is in the minimal degree of difficulty $\mindegdif$.
Intuitively, $\masspbofset(\mathcal{A})\massred\masspbofset(\mathcal{B})$ means
that there exists a Turing machine that computes an enumeration of $\mathcal{A}$
when given an enumeration of $\mathcal{B}$.

Given an alphabet $\alphab$ and a dimension $d$, one can fix a recursive
bijection between $\N$ and the set of finite patterns of dimension $d$ over the
alphabet $\alphab$.
Therefore we can represent these forbidden patterns as a $w\in\N^\N$:
$\val{w}{0}=\card{\alphab}$, $\val{w}{1}=d$ and
for $n\geq{}2$, $\val{w}{n}$ is the code of a forbidden pattern.
Effective subshifts can thus be expressed in terms of recursive operators: To an
integer $n$ associate the effective subshift $\effofcode{n}$ whose set of
forbidden patterns is coded in that way by
$w=\recop{n}(0^\infty)\in\N^\N$.

A set of effective subshifts $C$ is therefore coded by $\mathcal{C}$, a
subset of $\N$, and has its associated mass problem $\masspbofset(\mathcal{C})$.
Again, there exists a recursively enumerable set $K\subseteq\N$ such that
$C=\left\{\effofcode{n}| n\in{}K\right\}$ if and only if
$\masspbofset(\mathcal{C})$ is in $\mindegdif$.
Since any subshift can be written as the intersection of SFTs, we can also code
an arbitrary subshift in the same way: A set of effective subshifts
$\mathcal{C}\subseteq\N$ is said to code the subshift $\sshift{X}$ if
$\sshift{X}=\bigcap_{n\in{}\mathcal{C}}\effofcode{n}$.
Again, to a subshift $\sshift{X}$, we can associate a canonical mass problem
$\masspbofsubshift(\sshift{X})$ containing all its possible codings:
$\masspbofsubshift(\sshift{X})=\left\{x\in\N^{\N}|
\bigcap_{n\in\N}\effofcode{\val{x}{n}}=\sshift{X}\right\}
$.
As always, $\sshift{X}$ is an effective subshift if and only if there exists an
$n$ such that $\sshift{X}=\effofcode{n}$, \ie if and only if
$n^{\infty}\in\masspbofsubshift(\sshift{X})$, that is, if and only if
$\masspbofsubshift(\sshift{X})$ is in 
$\mindegdif$.

\subsection{Universality}

Since recursive operators may be too powerful for a sane definition of
universality, we may want to restrict them by considering a set of allowed
operations.
A set of integers $\mathcal{R}$ defines a \emph{reduction class} $R$:
$R=\left\{\recop{n}| n\in\mathcal{R}\right\}$.
Remark that, \eg the set of all factor maps forms a recursively enumerable
reduction class since they correspond to sliding block codes~\cite{Hedlund69}
and such codes can be recursively enumerated.
The same applies to the set of all sofic projective subdynamics.
We can now define what we call \emph{universality}:

\begin{definition}
        \label{defi:univ}
Let $\mathcal{R}$ be a reduction class, $\mathcal{C}$ and $\mathcal{G}$ sets of
effective subshifts and $\mathcal{B}$ a subset of $\N^2$, we say that a subshift
$\sshift{X}$ is \emph{\univclas} for $\mathcal{C}$ when:
\[
 i\in{}\mathcal{C} \Leftrightarrow\left\{
 \begin{array}{l}
         i\in{}\mathcal{G} \textrm{ and}\\
 \exists n\in \mathcal{R}, (i,b)\in\mathcal{B},
 \recop{n}(\sshift{X})_{[-b;b]^d} \textrm{ does not contain any} \\ \hfill\textrm{of the
 first }b\textrm{ forbidden patterns in }
 \effofcode{i}
 \end{array}\right.
\]

Similarly, we say that $\sshift{X}$ is \emph{\univclasinf} for $\mathcal{C}$ when
the following holds:
\[
 i\in{}\mathcal{C} \Leftrightarrow\left\{
 \begin{array}{l}
         i\in{}\mathcal{G} \textrm{ and}\\
 \exists n\in \mathcal{R},
 \recop{n}(\sshift{X})\subseteq\effofcode{i}
 \end{array}\right.
\]
\end{definition}

The set $\mathcal{B}$ shall be understood as the precision we require from the
simulation: The biggest $b$ is, closer the image of $\sshift{X}$ to elements of
$\effofcode{i}$ will be.
If $b$ could be chosen infinite then the image of $\sshift{X}$ would not contain
any $\effofcode{i}$-forbidden pattern.
This is why, in the second definition, we
replaced $\mathcal{B}$ by $\infty$ to emphasize that we require infinite
precision, \ie complete inclusion.
In some cases, with sufficiently big $b$'s, we can force the inclusion directly:
If $\mathcal{R}$ is the set of all factor maps, $\mathcal{G}$ the set of all
SFTs and for any $b_i$ such that $(i,b_i)$ is in $\mathcal{B}$, $b_i$ must be
big enough so that any forbidden pattern defining $\effofcode{i}$ is supported
on $[-b_i;b_i]^d$. This can actually be proved for any
subshift-preserving recursive operator:

\begin{lemma}
        \label{lemma:rgbunivinfuniv}
Given a set of SFTs $\mathcal{G}$, there exists
$\mathcal{B}\subseteq{}\N^2$ such that
$\masspbofset(\mathcal{B})\massred{}\masspbofset(\mathcal{G})$ and
for any subshift $\sshift{X}$ and any reduction class $\mathcal{R}$ of
$\sshift{X}$-subshift-preserving recursive operators, 
for any set of SFTs $\mathcal{C}$, $\sshift{X}$ is \univclas{} for
$\mathcal{C}$ if and only if  $\sshift{X}$ is \univclasinf{} for $\mathcal{C}$.
\end{lemma}

\begin{proof}
        We are given an enumeration of $i\in\mathcal{G}$ and want to compute an
        enumeration of the set $\mathcal{B}$.
        Let $b_i$ be a big enough integer such that any forbidden pattern of
        $\effofcode{i}$ is supported on $[-b_i;b_i]^d$, where $\effofcode{i}$ is
        an SFT of dimension $d$.
        When $i$ is given, this $b_i$ can be easily computed, so we define:
        $\mathcal{B}=\left\{(i,b_i), i\in\mathcal{G}\right\}$ and the
        relation $\masspbofset(\mathcal{B})\massred{}\masspbofset(\mathcal{G})$
        holds.

        Suppose $\sshift{X}$ is \univclas{} for $\mathcal{C}$.
        Let $(i,b_i)\in\mathcal{B}$ and
        $y\in\recop{n}(\sshift{X})$.
        We supposed $\recop{n}$ to be $\sshift{X}$-subshift-preserving:  For
        any $v\in\Z^d$, there exists $x_v\in\sshift{X}$ such that
        $\recop{n}(x_v)=\sigma_v(y)$.
        Since all the forbidden patterns of
        $\effofcode{i}$ are supported on $[-b_i;b_i]^d$,
        $y$ does not contain any pattern
        forbidden in $\effofcode{i}$. Since $y$ is arbitrary, we conclude that
        $\recop{n}(\sshift{X})\subseteq\effofcode{i}$, thus 
        $\sshift{X}$ is \univclasinf{} for $\mathcal{C}$.
        The other direction is trivial.
\end{proof}

This definition of universality is weak:
$\recop{n}(\sshift{X})\subseteq\effofcode{i}$ implies that any
subshift $\effofcode{i'}\in{}G$ such that $\effofcode{i}\subseteq\effofcode{i'}$
also verifies $\recop{n}(\sshift{X})\subseteq\effofcode{i'}$.
The artifact introduced by considering universality relative to a set of
subshifts $\mathcal{G}$ is what allows to avoid trivial necessary
conditions on the set $\mathcal{C}$ for the existence of an universal subshift.
The artifact introduced by the set $\mathcal{B}$ in the definition of
\univclas{} is there to make proofs simpler: In fact, we are interested in
subshift-preserving reduction classes and with Lemma~\ref{lemma:rgbunivinfuniv}
being \univclas{} is equivalent to \univclasinf{} for a carefully chosen
$\mathcal{B}$.

A stronger notion of universality would be to ask $\sshift{X}$ to simulate
every point of $\effofcode{i}$ by asking $\recop{n}(\sshift{X})=\effofcode{i}$
in the definition of universality, instead of only some of them
(what we get with
$\recop{n}(\sshift{X})\subseteq\effofcode{i}$).
We will see in Section~\ref{sec:cons} that even this weak notion of universality
imposes some conditions on the universal subshifts.
In Section~\ref{sec:constr} we will show how to construct universal subshifts for
this stronger notion of universality ($\recop{n}(\sshift{X})=\effofcode{i}$) as soon
as the necessary conditions developed in Section~\ref{sec:cons} are met.
Therefore, considering this weaker definition of universality only makes our
results stronger.

\section{Consequences of the existence of universal subshifts}

\label{sec:cons}

M. Hochman proved that there cannot exist any universal effective subshift (in
any dimension) for the set of non-empty SFTs of dimension at least $2$ by
proving that this would give an algorithm enumerating all the non-empty
$2-$dimensional SFTs~\cite{hochmanuniv}. This would contradict Berger's
theorem~\cite{bergerthesis} stating that it is undecidable, given a finite set
of patterns $\ForbPat$, to know if $\sshift{X}_{\ForbPat}$ is empty in
dimension $d\geq 2$ because the set of empty SFTs is recursively
enumerable~\cite{wan}.
Hochman's definition of universality is, in our context, considering the
reduction class of sofic projective subdynamics. With a very similar proof, we
generalize his result to our notion of universality and give some easy
consequences that were not obvious from the original result.


As we are dealing with effective continuous functions,
we can expect to
obtain classical topological results in an effective way; the one that we will
use later is that a continuous function over a compact set is uniformly
continuous, and that ``uniform continuity can be computed'':

\begin{lemma}[{Straightforward extension of \cite[Chapter~II.3]{odifreddi1},\cite{uspenski55,nerode57}}]
\label{lemma:compuccomp}
There exists a Turing machine with oracle such that when given as oracle an
element of $\masspbofsubshift(\sshift{X})$, where $\sshift{X}$ is a subshift of
$\myspace$,  and as input $n$ and $r$ computes, when
$\recop{n}$ is defined on $\sshift{X}$,
$l=\ucfunc(n,r,\sshift{X})$ such that for any
$\conf{x},\conf{y}\in\sshift{X}$, we have:
\[
        \conf{x}_{|[-l;l]^d}=\conf{y}_{|[-l;l]^d}\Rightarrow
\recop{n}(\conf{x})_{|[0;r]}=\recop{n}(\conf{y})_{|[0;r]}
\]
\end{lemma}

In simpler words: if two configurations on a finite alphabet match on a
``sufficiently big'' domain around the origin then their images have a long
common prefix; and this ``sufficiently big'' can be computed.
Note that this proof is classical, as it follows trivially from the proof of
classical results, but we have not been able to find a reference to use for
our case.
The same proof applies to compact subsets of $\N^{\N}$ but since we have the
formalism for dealing with subshifts, we only prove it in this case.

\begin{proof}
        Since $\recop{n}$ is continuous over a compact set $\sshift{X}$, it is
uniformly continuous: $\forall r\in\N, \exists l_r\in\N, \forall
\conf{x},\conf{y}\in{}\sshift{X}, \conf{x}_{|[-l_r;;l_r]^d}=\conf{y}_{|[-l_r;l_r]^d}\Rightarrow
\recop{n}(\conf{x})_{|[0;r]}=\recop{n}(\conf{y})_{|[0;r]}$.
This $l_r$, or any larger integer, is the $\ucfunc(n,r,\sshift{X})$ we are
looking for.

For an integer $i$ growing from $0$ to $\infty$, for any pattern $\motif{P}$ of
$\myspace$ defined on $[-i;i]^d$, check that $\motif{P}$ does not contain any of
the first $i$ forbidden patterns of $\sshift{X}$ given by 
the oracle $\masspbofsubshift(\sshift{X})$. 
If it does not contain such a pattern then
we simulate $\machine{\recop{n}}$ given as oracle
$\motif{P}$ on every input in $\left\{0,\ldots,r\right\}$.
In this simulation, we distinguish two cases:
\begin{enumerate}
\item If for at least one $\motif{P}$ and one input, before halting, the
        machine asks for the value of its oracle $\motif{P}$ outside
        of the domain of $\motif{P}$ ($[-i;i]^d$) then continue with the next $i$.
\item Otherwise define $\ucfunc(n,r,\sshift{X})=i$.
\label{item:uceffarret}
\end{enumerate}

By construction, $\ucfunc(n,r,\sshift{X})$
matches the conclusions of our lemma since the machine never reads anything
outside of $[-\ucfunc(n,r,\sshift{X});\ucfunc(n,r,\sshift{X})]^d$ and thus
values can be changed outside of this domain without altering the image on its
first $r$ values.
However, we have to prove that case~\ref{item:uceffarret} always
happens for a given $n,r$ and $\sshift{X}$ so that $\ucfunc$ is well defined.

Suppose case~\ref{item:uceffarret} does not happen: For any integer $i$, there
exists a pattern $\motif{P}_i$ with domain $[-i;i]^d$ that does not contain any of the
first $i$ forbidden patterns of $\sshift{X}$,
and an input
$p_i\in\left\{0,\ldots,r\right\}$ such that $\machine{\recop{n}}$ asks for a
value of its oracle outside of $[-i;i]^d$ with input $p_i$.
Since the input $p_i$ can take only finitely many different values, this implies:
There exists an input $p$ such that for infinitely many $i$'s, there
exists a $\motif{P}_i$ defined over $[-i;i]^d$ such that
$\machine{\recop{n}}$ asks for a value of its oracle outside of the domain of
$\motif{P}_i$ with input $p$.
Let $v_1\in\N$ be the first value asked by the Turing machine. Since we still
have infinitely many patterns $\motif{P}_i$, infinitely many of them have the
same value at $v_1$, thus, without loss of generality we assume that all of the
$\motif{P}_i$'s have the same value at $v_1$. Since the Turing machine is
deterministic and has a fixed input, the second value it will ask is the same
for every $\motif{P}_i$: $v_2$. Again, w.l.o.g., we may assume that all the
$\motif{P}_i$'s have the same value at $v_2$.
By iterating this process we obtain a configuration $\conf{c}\in\myspace$ such
that the Turing machine asks for the value $v_1$, then $v_2$, etc. when given
$\conf{c}$ as oracle. Since this machine asks for infinitely many values of its
oracle, it never halts. However, by its construction, $\conf{c}$ does not
contain any forbidden pattern for $\sshift{X}$, hence $\conf{c}\in\sshift{X}$
and $\recop{n}$ is not defined on all $\sshift{X}$.
\end{proof}

Note that in the proof of Lemma~\ref{lemma:compuccomp}, if $\recop{n}$ is not
defined on all $\sshift{X}$ then the algorithm may not halt, so that
$\ucfunc(n,r,\sshift{X})$ is only defined for $n$'s such that $\recop{n}$ is
defined on $\sshift{X}$.

We now have all the tools to generalize M. Hochman's proof~\cite{hochmanuniv}:

\begin{theorem}
\label{thm:univgtclass}
If a subshift $\sshift{X}$ is \univclas{} for a set of effective subshifts
$\mathcal{C}$ then:
\[
\masspbofset(\mathcal{C})\massred
\masspbofsubshift(\sshift{X})
\medvedevsup\masspbofset(\mathcal{R})
\medvedevsup\masspbofset(\mathcal{G})
\medvedevsup\masspbofset(\mathcal{B})
\]
\end{theorem}

\begin{proof}
We are given as oracle
$x\in
\masspbofsubshift(\sshift{X})
\medvedevsup\masspbofset(\mathcal{R})
\medvedevsup\masspbofset(\mathcal{G})
\medvedevsup\masspbofset(\mathcal{B})\subseteq\N^{\N}$;
by decoding $x$, we may assume that we are given $\alphab$ and $d$ such that $\sshift{X}$ is a subshift
of $\alphab^{\Z^d}$ defined by a family of forbidden patterns
$(\motif{F}_i)_{i\in\N}$, that we can enumerate given $x$.

To prove our theorem, we need to prove that we can enumerate all the subshifts
of $\mathcal{C}$ given such an oracle, which will give us an element of
$\masspbofset(\mathcal{C})$.
The following algorithm is an adaptation to our context of M. Hochman's
one~\cite[Section~$5$]{hochmanuniv}:

\begin{itemize}
        \item Enumerate all the integers $i\in\mathcal{G}$ and decode
                $k_i$, $d_i$ and $(\ForbPat(i)_n)_{n\in\N}$ which are, respectively, the size of the alphabet,
the dimension and the forbidden patterns defining
$\effofcode{i}$.
\item Enumerate all the integers $n$ in $\mathcal{R}$ and all the integers $b$ such
that $(i,b)\in\mathcal{B}$.
\item By hypothesis, $\recop{n}$ is defined on $\sshift{X}$ since it is
        $\sshift{X}$-subshift-preserving: Compute $l=\ucfunc(n,b,\sshift{X})$
        from Lemma~\ref{lemma:compuccomp}.
\begin{itemize}
\item For every $j>l$,  compute all the patterns defined on
$[-j;j]^d$ that do not contain any of the $\motif{F}_1,\ldots,\motif{F}_j$ and
denote them by $\motif{P}_1,\ldots\motif{P}_h$.
\item If all the patterns
        $\recop{n}(\motif{P}_1)_{|[-b;b]^{d}},\ldots,\recop{n}(\motif{P}_h)_{|[-b;b]^{d}}$
        do not contain any of the first $b$ forbidden patterns in
        $\effofcode{i}$ then claim that $\effofcode{i}\in{}\mathcal{C}$:
        Enumerate it.
\end{itemize}
\end{itemize}

By definition, this algorithm enumerates only elements of $\mathcal{C}$.
It remains to prove that we actually enumerate the whole set $\mathcal{C}$.
Suppose it is not the case:
If there exists $n\in\mathcal{R}$, $i\in \mathcal{C}$, $(i,b)\in\mathcal{B}$
such that $\recop{n}(\sshift{X})_{|[-b;b]^d}$ does not contain any of the first
$b$ forbidden patterns in $\effofcode{i}$ and the algorithm does not
enumerate $i$: For arbitrary large $j$, there is a pattern $\motif{P}_k$ not
containing any of the $\motif{F}_1,\ldots,\motif{F}_j$
defined on $[-j;j]^d$ such that
$\recop{n}(\motif{P}_k)_{|[-b;b]^{d}}$ contains one of the first $b$ forbidden
patterns in $\effofcode{i}$. By compactness, there would exist $\conf{x}\in
\sshift{X}$ such that $\recop{n}(x)_{|[-b;b]^{d}}$ contains one of the first $b$
forbidden patterns in $\effofcode{i}$.
\end{proof}

Combining it with Lemma~\ref{lemma:rgbunivinfuniv},
we get a result for \univclasinf{} subshifts:

\begin{theorem}
\label{thm:univgtclasssft}
For any set of SFTs $\mathcal{G}$, any subshift $\sshift{X}$ and any reduction
class of $\sshift{X}$-subshift-preserving recursive operators $\mathcal{R}$,
if
$\sshift{X}$ is \univclasinf{} for $\mathcal{C}$ then:
\[
\masspbofset(\mathcal{C})\massred
\masspbofsubshift(\sshift{X})
\medvedevsup\masspbofset(\mathcal{R})
\medvedevsup\masspbofset(\mathcal{G})
\]
\end{theorem}

\begin{proof}
 Find $\mathcal{B}$ from Lemma~\ref{lemma:rgbunivinfuniv};
$\masspbofset(\mathcal{B})\massred{}\masspbofset(\mathcal{G})$
and 
if $\sshift{X}$ is \univclas{} for
$\mathcal{C}$ then $\sshift{X}$ is \univclasinf{} for $\mathcal{C}$.
From Theorem~\ref{thm:univgtclass}:
$\masspbofset(\mathcal{C})\massred
\masspbofsubshift(\sshift{X})
\medvedevsup\masspbofset(\mathcal{R})
\medvedevsup\masspbofset(\mathcal{G})
\medvedevsup\masspbofset(\mathcal{B})$ but then since
$\masspbofset(\mathcal{B})\massred{}\masspbofset(\mathcal{G})$ we get our
result.
\end{proof}

In our formalism, M. Hochman considers $\mathcal{R}$ to be a recursively
enumerable set (the set of all sofic projective subdynamics maps) and
$\mathcal{G}$ to be the set of all (possibly empty) SFTs, which is also
recursively enumerable, the theorem can thus be reformulated as
$\masspbofset(\mathcal{C})\massred\masspbofsubshift(\sshift{X})$ with these
conditions and obtain Hochman's result~\cite{hochmanuniv}.
By considering $\mathcal{R}$ to be the set of recursive operators corresponding
to the $\preceq$ relation from \cite{gtlbunivac} we can also recover
\cite[Theorem~$1$]{gtlbunivac}.

Note that the condition on $\mathcal{R}$ is important:
Consider an SFT $\sshift{X}$ that is maximal in the lattice of
Medvedev degrees of $2-$dimensional SFTs which we know exists by Simpson's
work~\cite{simpsonmedv}. We also know that for any non-empty SFT $\sshift{S}$
there exists a recursive operator $\Psi$ such that
$\Psi(\sshift{X})=\sshift{S}$~\cite[Lemma~$3.16$, point~$1$]{simpsonprod}.
If we take $\mathcal{R}$ to be the reduction class of all such $\Psi$'s (they
all are $\sshift{X}$-subshift-preserving by definition),
$G$ to be the set of all SFTs, we get
$\masspbofsubshift(\sshift{X})\medvedevsup\masspbofset(G)\in\mindegdif$ but the
set $C$ of non-empty SFTs is not in $\mindegdif$ so that
$\masspbofset(\mathcal{R})\not\in\mindegdif$ by
Theorem~\ref{thm:univgtclasssft}.

\begin{corollary}
\label{cor:noequalentropuniv}
Given any (non-empty) SFT $\sshift{S}$, there exists no effective subshift that
factors (nor projects by sofic projective subdynamics) onto all the non-empty SFTs of
entropy equal to the entropy of $\sshift{S}$ in dimension greater than one.
\end{corollary}

\begin{proof}
By \cite[Corollary~$3.3$]{hochman2007cem}, let $M$ be a Turing machine
computing a right approximation of the entropy of $\sshift{S}$: $M$
enumerates a sequence of rational numbers $(h_n)_{n\in\N}$ such that
$\lim_{n\to\infty}h_n=\entrop{\sshift{S}}$ and $\forall n\in\N,
h_n\geq{}h_{n+1}\geq\entrop{\sshift{S}}$.

Let $M_n$ be the Turing machine executing $M$ one step out of two and
the $n^{th}$ Turing machine 
the other step. $M_n$ halts whenever one of the two machines halts, that is, $M_n$
halts if and only if the $n^{th}$ Turing machine halts.
With this Turing machine $M_n$, construct an SFT $\sshift{S}_n$ as in the proof
of \cite[Theorem~$1.1$]{hochman2007cem}, with a slight modification: If $M_n$
halts then we force a state that is incompatible with every other state so that
$\sshift{S}_n$ is empty. Thus, by \cite[Theorem~$1.1$]{hochman2007cem},
$\sshift{S}_n$ is empty if and only if $M_n$ halts and
$\entrop{\sshift{S}_n}=\entrop{\sshift{S}}$ otherwise.

Let $\mathcal{G}$ be the mass problem associated to the set
$\left\{\sshift{S}_n| n\in\N\right\}$.
Let $\mathcal{R}$ be the recursively enumerable reduction class of all factor
maps.
If we suppose that there exists an effective subshift factoring on all the SFTs
of entropy equal to the entropy of $\sshift{S}$, by
Theorem~\ref{thm:univgtclass} we get that the set $\left\{\sshift{S}_n|
\sshift{S}_n\neq\emptyset\right\}$ is recursively enumerable since
$\masspbofsubshift(\sshift{X})$, $\masspbofset(G)$ and
$\masspbofset(\mathcal{R})$ all belong to the minimal degree of difficulty
$\mindegdif$.
This would mean that the set of Turing machines that do not halt is recursively
enumerable, which is known not to be.
\end{proof}

\begin{corollary}
No effective subshift can factor onto all the zero-entropy SFTs.
\end{corollary}
\begin{proof}
Apply Corollary~\ref{cor:noequalentropuniv} with $\sshift{S}$ being the SFT
consisting of only one uniform configuration (this SFT has zero-entropy).
\end{proof}

We also note that an extended version of \cite{pierrebabispremin} (personal
communication) should prove that minimal non-empty SFTs are \emph{not}
recursively enumerable. Together with Theorem~\ref{thm:univgtclass}, this would
answer the question of M. Hochman whether or not there exists universal
effective subshifts for the set of non-empty minimal SFTs~\cite{hochmanuniv}:
Such universal subshifts would not exist.

\section{Construction of universal subshifts}

\label{sec:constr}

After proving necessary conditions for the existence of
universal subshifts, we construct universal subshifts that achieve
optimal conditions in the sense of Theorem~\ref{thm:univgtclass}.

\subsection{Toeplitz sequences as foundations}

The main ingredient of all the following results is summarized by
Theorem~\ref{thm:zeno}: There exists an effective way to code independently an
infinity of configurations into a single one and the decoding is effective
too.

\begin{theorem}
\label{thm:zeno}
There exists a finite alphabet $\alphab$ and an effective $1-$dimensional subshift
$\sshift{X}\subseteq\alphab^{\Z}$ such that each configuration of $\sshift{X}$
contains $\omega$ independent layers, in the following sense:
\begin{enumerate}
\item 
There exists a recursively
enumerable set of factor maps
$(\varphi_n:\alphab^{\Z}\to(\couleur{\left\{0,1\right\}^n})^{\Z})_{n\in\N}$ and a computable sequence
of integers $(m_n)_{n\in\N}$ defining finite index sofic projective subdynamics maps
$\layer{n}: \alphab^{\Z}  \to  (\couleur{\left\{0,1\right\}^n})^{\Z}$ such
that $\val{\layer{n}(c)}{i}=\val{\varphi_n(c)}{m_ni}$.
\label{coding:facteur}
\item Each $\layer{n}$ is a surjection from $\sshift{X}$ onto
$(\couleur{\left\{0,1\right\}^n})^{\Z}$
\label{coding:surj}
\item There exists an algorithm such that when given as input a (possibly
infinite) list of forbidden patterns $\ForbPat$, defining a $1-$dimensional
subshift $\sshift{X}_{\ForbPat}$, the algorithm computes an integer $n_i$ and
enumerates an infinite list of forbidden patterns $\mathcal{A}(\ForbPat)$
defining a $1-$dimensional subshift $\sshift{X}_{\mathcal{A}(\ForbPat)}$ such
that:
\[\layer{n_i}(\sshift{X}\cap\sshift{X}_{\mathcal{A}(\ForbPat)})=\sshift{X}_{\ForbPat}\]
and, for $j\neq{}n_i$, if $\sshift{X}_{\mathcal{A}(\ForbPat)}\neq\emptyset$:
\[\layer{j}(\sshift{X}\cap\sshift{X}_{\mathcal{A}(\ForbPat)})=\layer{j}(\sshift{X})\]
\label{coding:indep}
\end{enumerate}
\end{theorem}

The intuition behind Theorem~\ref{thm:zeno} is as follows:
$\sshift{X}$ forces each of its configurations to contain
an infinity of independent layers; we use $m_n$ cells in $\sshift{X}$ to code $n$
bits in the $n^{th}$ layer.
$\layer{n}$ is the decoding function: It reads the $n^{th}$ layer, gets the $n$
bits value and decodes it to get a configuration over
$(\couleur{\left\{0,1\right\}^n})^{\Z}$.
Point~\ref{coding:indep} states that we can forbid coded patterns in any layer
independently of the other layers.

\begin{proof}

Our proof is in four parts: First we describe how to code $\omega$ independent
layers and then we show how to code required information for our theorem, how to
decode it and finally why the layers are independent.

\textbf{Layers}:
Since each configuration contains an infinite number of cells, we use half of
our cells to code the first layer, then half of what is remaining for the second
layer, etc.
In fact, the base construction is a Toeplitz sequence constructed as in
\cite{toeplitz,toeplitzprop} in which each step of the construction is used to
obtain a new independent layer.

Fix an integer $k$ greater than $2$.
A meta-cell is a group of $k$ cells
(the \tikz[scale=.2]{\codingcell{0}{0};} cells in Figure~\ref{fig:level1}) bordered by left brackets
(\tikz[scale=.2]{\leftbrack{0}{0};}) on their left and right brackets
(\tikz[scale=.2]{\rightbrack{0}{0};}) on their right.
To the right of a left bracket (\tikz[scale=.2]{\leftbrack{0}{0};}) there are
exactly $k$ coding cells (\tikz[scale=.2]{\codingcell{0}{0};}) and then a right
bracket (\tikz[scale=.2]{\rightbrack{0}{0};}).
To the right of a right bracket (\tikz[scale=.2]{\rightbrack{0}{0};}) there are
exactly $k+2$ blank cells
(\tikz[scale=.2]{\draw (0,0) rectangle +(1,1);}) and then a left bracket
(\tikz[scale=.2]{\leftbrack{0}{0};}).
For
simplicity we use $k+2$ here because it is the size of a meta-cell.
This gives us a set of forbidden words for coding the first layer as depicted
in Figure~\ref{fig:level1}.

\begin{figure}[htb]
\begin{center}
\opt{dessinetikz}{\begin{tikzpicture}[scale=.25]

\draw (0,0) grid +(45,1);

\foreach \x in {0,12,24,36} {
\leftbrack{\x}{0};
}

\foreach \x in {5,17,29,41} {
\rightbrack{\x}{0};
}

\foreach \x in {1,13,25,37} {
\foreach \y in {0,1,2,3} {
\codingcell{\x+\y}{0};
}
}
\end{tikzpicture}}\opt{pdftikz}{\opt{bw}{\includegraphics{dessins/level1.pdf}}\opt{color}{\includegraphics{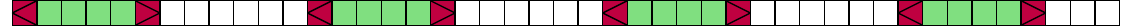}}}

\vspace{.3cm}

\begin{tabular}{ccl}
\tikz[scale=.25]{\leftbrack{0}{0};}&:& Left bracket\\
\tikz[scale=.25]{\rightbrack{0}{0};}&:& Right bracket\\
\tikz[scale=.25]{\codingcell{0}{0};}&:& Coding cell (described later)\\
\tikz[scale=.25]{\draw (0,0) rectangle +(1,1);}&:& Cell intentionally left blank\\
\end{tabular}

\end{center}
\caption{First layer (with $k=4$).}
\label{fig:level1}
\end{figure}

Now we will fill the blank cells (\tikz[scale=.2]{\draw (0,0) rectangle
+(1,1);}) in order to code the second layer.
First, remark that if we group the cells by arrays of $k+2$ cells we see that
half the cells are occupied by the first layer and the other half is left free
as depicted on Figure~\ref{fig:level1meta}.

\begin{figure}[htb]
\begin{center}
\opt{dessinetikz}{\begin{tikzpicture}[scale=.25]

\draw (0,0) grid +(45,1);

\foreach \x in {0,2,4,...,44} {
\obscell{\x}{0};
}

\end{tikzpicture}}\opt{pdftikz}{\opt{bw}{\includegraphics{dessins/level1meta.pdf}}\opt{color}{\includegraphics{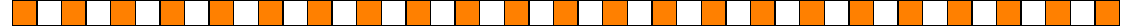}}}

\vspace{.3cm}

\begin{tabular}{ccl}
\tikz[scale=.25]{\obscell{0}{0};}&:& Cell occupied by the first layer\\
\tikz[scale=.25]{\draw (0,0) rectangle +(1,1);}&:& Cell left free by the first layer\\
\end{tabular}

\end{center}
\caption{Occupied and free cells after the first layer (here one cell represents
a block of $k+2$ real cells).}
\label{fig:level1meta}
\end{figure}

Now, we see the discrete line as divided in blocks of size $2(k+2)$, each of
these blocks having $k+2$ blank cells.
For the second layer, we use meta-cells consisting of $k+2$ identical real
cells. Those meta-cells will fill the cells left blank by the first layer.
After that, this is basically the same construction as for the first layer. 
A \{left bracket, right bracket, coding cell, blank cell\} for the second layer
consists of $k+2$ ``real'' \{left brackets, right brackets, coding cells, blank
cells\}.
The same rules as for the cells of the first layer apply to these meta-cells.
The global picture of a configuration with the first and second layers is
depicted in Figure~\ref{fig:level2}.

\begin{figure}[htb]
\begin{center}
\opt{dessinetikz}{\begin{tikzpicture}[scale=.12]

\draw (0,0) grid +(90,1);

\foreach \x in {0,2,4,...,88} {
\obscell{\x}{0};
}

\foreach \x in {1,25,49,73} {
\leftbrack{\x}{0};
}

\foreach \x in {11,35,59,83} {
\rightbrack{\x}{0};
}

\foreach \x in {3,27,51,75} {
\foreach \y in {0,2,4,6} {
\codingcell{\x+\y}{0};
}
}

\node[draw] (obs) at (1,-10) {
\tikz[scale=.25]{
\leftbrack{0}{0};
\foreach \x in {1,2,3,4} {
\codingcell{\x}{0};
}
\rightbrack{5}{0};
}
};

\draw (2.5,0.5) -- (obs);

\node[draw] (cod2) at (60,-10) {
\tikz[scale=.25]{
\foreach \x in {0,1,2,3,4,5} {
\codingcell{\x}{0};
}
}
};

\draw (53.5,0.5) -- (cod2);

\node[draw] (lb2) at (15,10) {
\tikz[scale=.25]{
\foreach \x in {0,1,2,3,4,5} {
\leftbrack{\x}{0};
}
}
};

\draw (25.5,0.5) -- (lb2);

\node[draw] (rb2) at (70,10) {
\tikz[scale=.25]{
\foreach \x in {0,1,2,3,4,5} {
\rightbrack{\x}{0};
}
}
};

\draw (83.5,0.5) -- (rb2);

\end{tikzpicture}}\opt{pdftikz}{\opt{bw}{\includegraphics{dessins/level2.pdf}}\opt{color}{\includegraphics{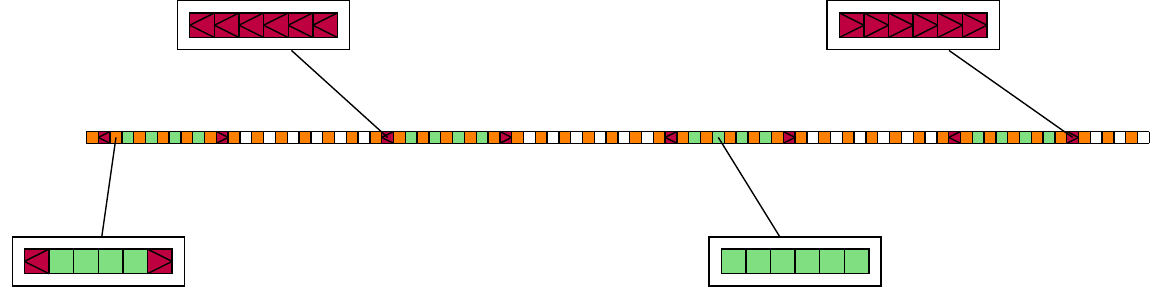}}}

\vspace{.3cm}

\begin{tabular}{ccl}
\tikz[scale=.25]{\obscell{0}{0};}&:& Cell occupied by the first layer\\
\tikz[scale=.25]{\leftbrack{0}{0};}&:& Left bracket\\
\tikz[scale=.25]{\rightbrack{0}{0};}&:& Right bracket\\
\tikz[scale=.25]{\codingcell{0}{0};}&:& Coding cell (described later)\\
\tikz[scale=.25]{\draw (0,0) rectangle +(1,1);}&:& Cell intentionally left blank\\
\end{tabular}

\end{center}
\caption{First and second layers.}
\label{fig:level2}
\end{figure}

Note that adding the second layer does not affect the first one if we modify a
bit our rules: A left bracket of the first layer is the only one with a
coding cell at its right and a right bracket has a coding cell on its left;
hence we can still force a given number of blank and coding cells
around the first layer brackets by considering the cells it expects to be blank
as blank.

The second layer, as the first one, occupies blocks of $k+2$ meta-cells and
leaves $k+2$ meta-blank cells between them. By dividing again the line in blocks
of $2(k+2)$ meta-cells, each block contains $k+2$ meta-blank cell; that is, each
block consists of $(2(k+2))^2$ real cells and contains $(k+2)^2$ real blank
cells.
The third layer is the same as the second one except we use meta-meta-cells of
meta-cells, and so on for each layer.

We described a procedure to produce the forbidden patterns defining
a subshift, hence the obtained subshift is effective.
As a conclusion of this part of the proof, we have:
\begin{itemize}
\item $\omega$ layers in each configuration of this effective subshift;
\item cells belonging to a given layer can be recognized locally: The
window needed depends only on the layer of cells we are looking at;
\item each meta-coding cell of the $n^{th}$ layer is an array of
$k(k+2)^{n-1}$ real coding cells
\end{itemize}

\textbf{Coding}:
The real coding cells (the \tikz[scale=.2]{\codingcell{0}{0};} cells in
Figure~\ref{fig:level1}), in fact,
are two: $\tikz[scale=.2]{\codingcell{0}{0};}_0$ and
$\tikz[scale=.2]{\codingcell{0}{0};}_1$ and thus code one bit of
information.
Remark that a meta-coding cell from the $n^{th}$ layer contains $k(k+2)^{n-1}$
real coding cells.
The leftmost cell of meta-coding cells in the $n^{th}$ layer are separated by
$2^{n}(k+2)^{n}$ cells, we therefore define $m_n = 2^{n}(k+2)^{n}$ in
point~\ref{coding:facteur}.

\textbf{Decoding}:
In this paragraph, we define the factor maps $\varphi_n$ for proving our
point~\ref{coding:facteur}.
We define $\varphi_n$ as a block map; $\varphi_n$ has a radius of $m_n/2$.
$\varphi_n$, when applied to a configuration of our subshift, at every cell,
``sees'' at least one entire meta-coding cell from the $n^{th}$ layer: choose,
\eg the leftmost one.
A meta-coding cell from the $n^{th}$ layer codes $k(k+2)^{n-1}$ bits of
information.
We want $\varphi_n$ to take its values in $\left\{0,1\right\}^n$ in
point~\ref{coding:facteur}; since $k>2$, $k(k+2)^{n-1} > n$, and there exists a
surjection from $\left\{0,1\right\}^{k(k+2)^{n-1}}$
onto $\left\{0,1\right\}^n$ which will be the value of $\varphi_n$ at this
position. Note that we lose a lot of information with this
factor map but it does not matter since $n$ can be taken as large as we want.

The construction of $\varphi_n$ and $m_n$ is effective, which proves
point~\ref{coding:facteur}.
Since, with these $\varphi_n$ and $m_n$, $\layer{n}$ is a
surjection from $\sshift{X}$ onto $(\couleur{\left\{0,1\right\}^n})^{\Z}$,
this proves point~\ref{coding:surj}.

\textbf{Independence}:
Given an enumeration $\ForbPat$ of patterns defining a $1-$dimensional subshift,
defined on an alphabet of size smaller than $2^n$, for a given integer $n$, it
is clear that this alphabet can be sent into $\left\{0,1\right\}^n$ and decoded
back.
$n$, or any greater integer, is thus the $n_i$ we need for
point~\ref{coding:indep}.
Given an enumeration of $\ForbPat$, we can effectively transform these patterns
so that they correspond to their codings on the $n^{th}$ layer; these
transformed patterns are $\mathcal{A}(\ForbPat)$.
Since $\layer{n}$ depends only on the value of the configuration on its
$n^{th}$ layer, this completes the proof of point~\ref{coding:indep} and of
the whole theorem.
\end{proof}

We can reformulate Point~\ref{coding:facteur} as follows:
$m_n\Z$ is a subgroup of $\Z$ and define
$h_n(i)=m_ni:\Z\to{}m_n\Z$. $L_n(\sshift{X})$ is simply
$\projsub{\sshift{X}}{h_n}{\varphi_n}$, the $(h_n,\varphi{}_n)$
finite index sofic projective subdynamics of $\sshift{X}$.
Also, from the definition of the $L_n$'s, these define a recursively
enumerable reduction class $\mathcal{R}=\left\{r_n| n\in\N\right\}\subseteq\N$
of sofic projective subdynamics maps.

\subsection{Universality with subshift-preserving recursive operators}

We start by an almost straightforward application of Theorem~\ref{thm:zeno} to
construct a one-dimensional universal subshift with sofic projective
subdynamics as soon as the necessary conditions on the dimension and the
computability conditions imposed by Theorem~\ref{thm:univgtclass} are fulfilled.

\begin{theorem}
\label{thm:existsunivdim1}
For any set of non-empty $1-$dimensional effective subshifts $C$, there exists
a $1-$dimensional subshift $\sshift{X}$ 
such that $
        \masspbofsubshift(\sshift{X})\massred{}\masspbofset(\mathcal{C})$ and
there exists a recursively enumerable reduction class $\mathcal{R}$
such that $C=\left\{\recop{n}(\sshift{X})| n\in\mathcal{R}\right\}$
and all the $\recop{n}$'s are (finite index) sofic
projective subdynamics maps.
\end{theorem}

By Theorem~\ref{thm:zeno}, we know how to code an infinite number of
configurations in any configuration of an effective subshift and how to keep the
coding and decoding effective.
The idea of this proof is that every element of $\sshift{X}$ must in some
sense contain an element of any subshift of $C$, that is where the $\omega$
layers will play their key role.

\begin{proof}
First, let us fix some notations: $\mathcal{C}$ codes a set of non-empty
effective subshifts $C$: $\mathcal{C}=\left\{c_i| i\in\N\right\}$
$\left(C=\left\{\effofcode{c}|c\in\mathcal{C}\right\}\right)$.
$\effofcode{c_i}$ is a subshift of $\alphab_i^{\Z}$ defined by the (recursively
enumerable) set of forbidden patterns $\ForbPat_i$.
Since, by definition, $\alphab_i$ and $\ForbPat_i$ can be recursively enumerated
from $c_i$, Theorem~\ref{thm:zeno} is almost everything we need.
The sofic projective subdynamics and the recursively enumerable reduction class
$\mathcal{R}$ are given by point~\ref{coding:facteur} in
Theorem~\ref{thm:zeno}.

We enumerate the forbidden patterns defining the subshift $\sshift{X}$, knowing
the $c_i$'s as follows:
Enumerate the patterns defining the effective
subshift $\sshift{X}$ from Theorem~\ref{thm:zeno}.
Assume that we have coded up to $c_{i-1}$ and filled up to the $(n-1)^{th}$ layer
of our subshift. We are trying to code $c_i$ into the $n^{th}$ layer.
If $n$ bits are not enough to code $\alphab_{i}$ then code $c_{i-1}$ in the
$n^{th}$ layer and continue with the next layer.
This does not work for $c_1$
but this is only a detail: We may use a special code telling that this layer
codes nothing and still continue with the next layer, or, alternatively, choose
the initial length $k$ of the first layer in Theorem~\ref{thm:zeno} to be big
enough to contain it.

Now assume now that $n$ bits are enough to code $\alphab_{i}$.
Enumerate the patterns of $\ForbPat_{i}$ and transform them as in
point~\ref{coding:indep} of Theorem~\ref{thm:zeno} to get
$\mathcal{A}(\ForbPat_{i})$.
All these forbidden patterns are enumerated in parallel and define a
subshift $\sshift{X}$, which is non-empty since we assumed $C$ to contain only
non-empty subshifts.
It is clear that
$\masspbofsubshift(\sshift{X})\massred{}\masspbofset(\mathcal{C})$ and
point~\ref{coding:indep} of Theorem~\ref{thm:zeno} ensures us that
$
C=
\left\{\effofcode{c}|c\in\mathcal{C}\right\}
=
\left\{\layer{n}(\sshift{X})|n\in\N\right\}
$. 
\end{proof}

We now construct universal subshifts where we do not impose
any restriction on the reduction class, besides being subshift-preserving
recursive operators.
This construction allows us to get a general result about universality, proving
that in the context of subshift-preserving recursive operators, dimension does
not matter: As long as the necessary conditions from
Theorem~\ref{thm:univgtclass} are fulfilled, there exists a $1$-dimensional
universal subshift!
This shall be understood as a converse to 
Theorem~\ref{thm:univgtclass}.

\begin{theorem}
\label{thm:univ1d}
There exists a recursively enumerable reduction class $\mathcal{R}$ such that
for any set of non-empty effective subshifts $C$, 
there exists a $1-$dimensional subshift $\sshift{X}$ such that $\mathcal{R}$ is
$\sshift{X}$-subshift-preserving,
$\masspbofsubshift(\sshift{X})\massred{}\masspbofset(\mathcal{C})$ and
$C=\left\{\recop{n}(\sshift{X})|n\in\mathcal{R}\right\}$.
\end{theorem}

\begin{proof}
There exists a recursive injection of
$\bigcup_{d\in\N,\alphab\textrm{ finite}}\myspace$ into
$\left\{0,1\right\}^{\N}$.
With the construction of Theorem~\ref{thm:zeno}, there exists an effective
$1-$dimensional subshift $\sshift{K}$ where $\left\{0,1\right\}^{\N}$ can be
recursively injected into $\sshift{K}$: given $x\in\left\{0,1\right\}^{\N}$,
code $x_i$ as a uniform configuration in the $i^{th}$ layer of the construction.
Therefore, there exists a recursive operator $\Psi$ such that for any
effective subshift $\sshift{S}$ of $\myspace$ (where $\alphab$ and $d$ are not
fixed \latin{a priori}), we can construct an effective $1-$dimensional subshift
$\sshift{K_S}$ such that $\Psi(\sshift{K_S})=\sshift{S}$.
$\Psi$ is therefore $\sshift{K_S}$-subshift-preserving.
There exists an algorithm transforming the forbidden patterns defining
$\sshift{S}$ into those defining $\sshift{K_S}$.
Denote by $\varphi$ the computable function from $\N$ to $\N$ such that for all
$n$, $\effofcode{\varphi(n)}=\sshift{K_{\effofcode{n}}}$.
Applying Theorem~\ref{thm:existsunivdim1} to $\varphi(C)$ allows us to
conclude since we can get back to $C$ \latin{via} $\Psi$ which is
$\sshift{K_S}$-subshift-preserving for all $\sshift{K_S}$.
\end{proof}

If we consider a recursively enumerable set of non-empty effective subshifts,
Theorem~\ref{thm:univ1d} gives an effective universal subshift of dimension~$1$.
We know that effective subshifts of dimension~$1$ can be ``implemented'' by 
$2-$dimensional SFTs~\cite{mathieunathaliesoficeff,drseffsofic} in the sense that
the $1-$dimensional effective subshift is the sofic projective subdynamics of
the $2-$dimensional SFT.
We can therefore obtain the following corollary to Theorem~\ref{thm:univ1d}:
\begin{corollary}
\label{cor:univsft2d}
There exists a recursively enumerable reduction class $\mathcal{R}$ such that
for any recursively enumerable set of non-empty effective subshifts $C$ there
exists a $2-$dimensional SFT $\sshift{S}$ 
such that $\mathcal{R}$ is
$\sshift{S}$-subshift-preserving and
$C=\left\{\recop{n}(\sshift{S})|n\in\mathcal{R}\right\}$.
\end{corollary}

For example, one may take $C$ to be the set of non-empty $2-$dimensional SFTs
with a periodic point, of non-empty $1-$dimensional SFTs, of SFTs of dimension
prime with $42$ with a uniform configuration, limit sets of cellular automata,
etc.

\subsection{Universality with sofic projective subdynamics}

\label{sec:factuniv}

In Corollary~\ref{cor:univsft2d}, one may want to replace ``subshift-preserving
recursively enumerable reduction class'' by simply ``sofic projective
subdynamics'': This restriction is usually imposed in 
order not to hide the ``simulation power'' of the subshift in the reduction
itself.
This is, of course, not possible if we allow any dimension since the dimension
of the sofic projective subdynamics of a subshift $\sshift{X}$ is never greater
than the dimension of $\sshift{X}$. Even if we only
consider subshifts of dimension $2$, the finite index sofic projective
subdynamics of an SFT would be a sofic subshift and there exist non-sofic
effective subshifts, so that Corollary~\ref{cor:univsft2d} cannot be true if we
do not impose the set $C$ of simulated subshifts to be sofic.
If, moreover, we impose $C$ to be a set of SFTs (or sofic shifts) then it is
possible and has already been done by Lafitte and
Weiss~\cite[Theorem~$10$]{mwjournal}: In our formalism, this theorem states
precisely that given any recursively enumerable set of non-empty
two-dimensional SFTs, there exists a two-dimensional SFT that admits all these
SFTs as finite index sofic projective subdynamics. Because of
Theorem~\ref{thm:univgtclass}, this result is optimal since there is no hope to
be universal for a set of subshifts that is not recursively enumerable and the
sofic projective subdynamics of a non-empty subshift is necessarily non-empty.

In this section, we give a converse to Theorem~\ref{thm:univgtclass} with sofic
projective subdynamics and, in some sense, extend the aforementioned result of
Lafitte and Weiss to any set of non-empty effective subshifts: As soon as the set of simulated effective subshifts is
of bounded dimension $d$, and the computability conditions of
Theorem~\ref{thm:univgtclass} are fulfilled, there exists a universal subshift
of dimension $d$ for this set of effective subshifts.

\begin{theorem}
\label{thm:existsunivdimqcq}
For any set of non-empty $d-$dimensional effective subshifts $C$ there exists
a $d-$dimensional subshift $\sshift{X}$ such that any subshift of $C$ is
realized as the finite index sofic projective subdynamics of $\sshift{X}$ and 
$\masspbofsubshift(\sshift{X})\massred{}\masspbofset(C)$.
\end{theorem}

\begin{proof}

We extend the construction from Theorem~\ref{thm:zeno} to any
dimension by extending the layered construction and the way we code the
effective subshifts.

\textbf{Extending the layered construction:}
Put the layered construction from Theorem~\ref{thm:zeno} onto one axis of
$\Z^d$.
On all the other axes, impose that the cells have a constant type: A 
\tikz[scale=.25]{\leftbrack{0}{0};} cell must be next to a
\tikz[scale=.25]{\leftbrack{0}{0};} cell, etc.
For example, in dimension~$2$, we would obtain configurations like the one
depicted in Figure~\ref{fig:level2dim2}.

\begin{figure}[htb]
\begin{center}
\opt{dessinetikz}{\begin{tikzpicture}[scale=.23]

\draw (0,0) grid +(50,10);

\foreach \z in {0,1,...,9 } {

\foreach \x in {0,2,4,...,48} {
\obscell{\x}{\z};
}

\foreach \x in {1,25,49} {
\leftbrack{\x}{\z};
}

\foreach \x in {11,35} {
\rightbrack{\x}{\z};
}

\foreach \x in {3,27} {
\foreach \y in {0,2,4,6} {
\codingcell{\x+\y}{\z};
}
}

}

\node[draw] (obs) at (4,-5) {
\tikz[scale=.25]{
\leftbrack{0}{0};
\foreach \x in {1,2,3,4} {
\codingcell{\x}{0};
}
\rightbrack{5}{0};
}
};

\draw (2.5,0.5) -- (obs);

\node[draw] (cod2) at (20,-5) {
\tikz[scale=.25]{
\foreach \x in {0,1,2,3,4,5} {
\codingcell{\x}{0};
}
}
};

\draw (29.5,0.5) -- (cod2);

\node[draw] (lb2) at (15,15) {
\tikz[scale=.25]{
\foreach \x in {0,1,2,3,4,5} {
\leftbrack{\x}{0};
}
}
};

\draw (25.5,9.5) -- (lb2);

\node[draw] (rb2) at (40,15) {
\tikz[scale=.25]{
\foreach \x in {0,1,2,3,4,5} {
\rightbrack{\x}{0};
}
}
};

\draw (49.5,9.5) -- (rb2);

\end{tikzpicture}}\opt{pdftikz}{\opt{bw}{\includegraphics{dessins/level2dim2.pdf}}\opt{color}{\includegraphics{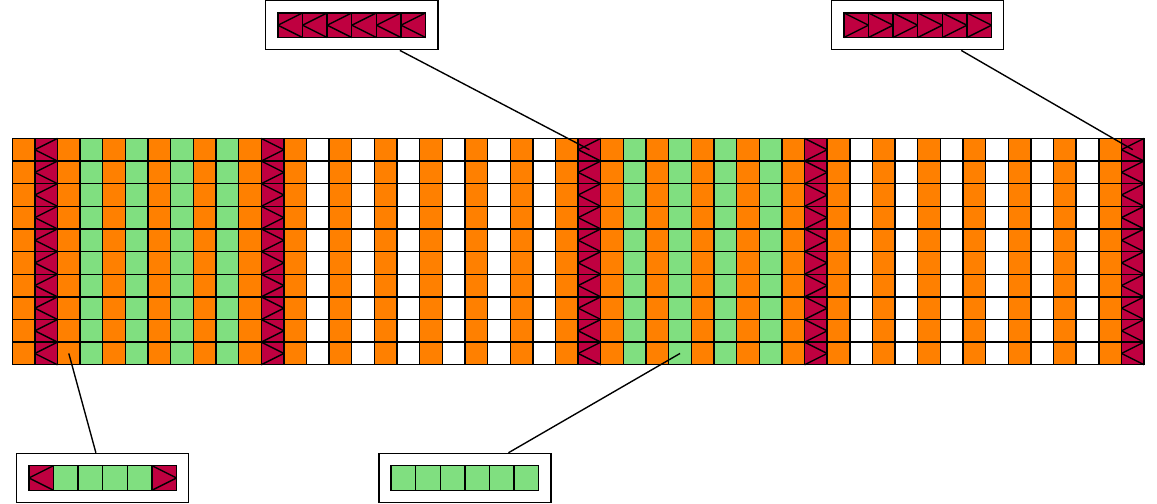}}}

\vspace{.3cm}

\begin{tabular}{ccl}
\tikz[scale=.25]{\obscell{0}{0};}&:& Cell occupied by the first layer\\
\tikz[scale=.25]{\leftbrack{0}{0};}&:& Left bracket\\
\tikz[scale=.25]{\rightbrack{0}{0};}&:& Right bracket\\
\tikz[scale=.25]{\codingcell{0}{0};}&:& Coding cell\\
\tikz[scale=.25]{\draw (0,0) rectangle +(1,1);}&:& Cell intentionally left blank\\
\end{tabular}

\end{center}
\caption{First and second layers in dimension~$2$.}
\label{fig:level2dim2}
\end{figure}

\textbf{Extending the coding and decoding of the subshifts:}
The subgroup of $\Z^d$ we consider is $m_n\Z\times\Z^{d-1}$ where
$m_n$ is from Point~\ref{coding:facteur} in Theorem~\ref{thm:zeno} and
the sofic projective subdynamics consists in applying $L_n$ independently on
every row of the configuration.
The forbidden patterns in the universal subshift are now those that are mapped
by this transformation to a forbidden pattern in $\sshift{X}_i\in{}C$ when
trying to code $\sshift{X}_i$ on the $n^{th}$ layer.
With these adaptations, the same proof as for Theorem~\ref{thm:existsunivdim1}
works.
\end{proof}

We are now able to construct universal subshifts for sets of subshifts that
all have the same dimension; since we can easily go down in dimension, we state
the result here:

\begin{theorem}
\label{thm:existsunivdimqcqdown}
For any set $C$ of non-empty effective subshifts of bounded dimension $d$,
there
exists $\sshift{X}$, a subshift of dimension $d$, such that any subshift of $C$ is
realized as the projective subdynamics of $\sshift{X}$ and
$\masspbofsubshift(\sshift{X})\massred{}\masspbofset(C)$.
\end{theorem}

\begin{proof}
The coding bits in the construction of Theorem~\ref{thm:zeno} take
their values in $\left\{0,1\right\}$.
Since the dimension is bounded and known we can make them take their values in
$\left\{0,1\right\}\times\left\{1,\ldots, d\right\}$; the second part will be
called the coded dimension. Coding bits that belong to the same layer must
have the same coded dimension.
If the coded dimension on a given layer is $k$, then the coding bits must be equal in
this layer on the $d-k+1$ last axes of $\Z^d$.
That way, we encode a $k-$dimensional subshift in a $d-$dimensional subshift in the
first $k$ dimensions of $\Z^d$ and can obtain the wanted projective subdynamics
by taking a subgroup of $\Z^d$ that ``forgets'' about the useless dimensions:
$m_n\Z\times\Z^{k-1}\times\left\{0\right\}^{d-k}$.
\end{proof}

\section{Conclusions}

In this paper we proved three main results: First, the Medvedev degree of 
the forbidden patterns defining a subshift cannot be lower than the degree of
the subshifts it simulates (Theorems~\ref{thm:univgtclass}
and \ref{thm:univgtclasssft}).
Second, that this necessary condition is optimal by giving a general
construction for obtaining a universal subshift given a set of effective
subshifts to simulate (Theorem~\ref{thm:univ1d} is the strongest result in
this sense).
Finally, by extending our construction of universal subshifts, we showed that
the simulation used can be restricted to sofic projective subdynamics by
obtaining the equivalent of Theorem~\ref{thm:univ1d} with
Theorem~\ref{thm:existsunivdimqcqdown} (the most general result of
section~\ref{sec:factuniv}).
Theorem~\ref{thm:existsunivdimqcqdown} is again optimal as the conditions on the
dimension are imposed by the nature of sofic projective subdynamics and the
conditions on Medvedev degree are imposed by
Theorem~\ref{thm:univgtclass}.

One may want to continue the weakening of the reduction class by replacing
``sofic projective subdynamics'' with ``factor'' in
Theorem~\ref{thm:existsunivdimqcq}.
Such a general statement cannot be true because the entropy of a factor of a
subshift is never greater that the entropy of the subshift itself.
Also, it seems difficult to use our construction from Theorem~\ref{thm:zeno}
since the ratio
of coding cells divided by the distance between them tends to zero as we
consider increasing layers.
Nevertheless, we can ask if bounding the entropy is again sufficient for the
existence of a universal subshift:
\begin{op}
Given a set of non-empty $d-$dimensional effective subshifts $C$ of entropy
bounded by $h$, does there exist a $d-$dimensional subshift $\sshift{X}$ such that
any subshift of $C$ is a factor of $\sshift{X}$ and
$\masspbofsubshift(\sshift{X})\massred{}\masspbofset(C)$?
Such that $\sshift{X}$ factors onto a subsystem of every subshift of $C$?
Can $\sshift{X}$ be of entropy $h$?
\end{op}

\bibliographystyle{amsplain}
\bibliography{biblio/biblio}

\end{document}